\documentclass[11pt]{article}
\usepackage[latin1]{inputenc}
\usepackage{amsmath}
\usepackage{mathtools}
\usepackage{comment}
\usepackage{upgreek}
\usepackage{multicol}
\usepackage{amsfonts}
\usepackage{mathrsfs}
\usepackage{bm}
\usepackage{gensymb}
\usepackage{wrapfig}
    \input{insbox}
\usepackage{enumerate}
\usepackage{enumitem}
\newcommand{\bfcalf}{\mbox{\boldmath ${\mathcal F}$}}
\newcommand{\bfcalg}{\mbox{\boldmath ${\mathcal G}$}}    
\usepackage{amssymb}
\usepackage{titlesec}
    \titleformat{\section}{\large\bfseries}{\thesection.}{0.5em}{}
    \titlespacing\section{0pt}{12pt plus 4pt minus 2pt}{6pt plus 2pt minus 2pt}
\usepackage{marvosym}
\usepackage{array}
\usepackage{chngcntr}
\usepackage{blindtext}
\usepackage[english]{babel}
\usepackage{amsthm}
    \theoremstyle{definition}
    \newtheorem{thm}{Theorem}[section]
    
    \newtheorem{lemma}[thm]{Lemma}
    \newtheorem{remark}[thm]{Remark}
    
    \newtheorem{definition}[thm]{Definition}
\usepackage{hyperref} 
    \hypersetup{
    colorlinks=true,
    linkcolor=blue,
    filecolor=blue,      
    urlcolor=blue,
    citecolor=blue,
    }
\usepackage{graphicx}
\usepackage[labelfont=bf]{caption}
\usepackage{color}
\usepackage{authblk}
\usepackage[titletoc,title]{appendix}
\usepackage{cite}
\usepackage[top=1in, bottom=1in, left=1in, right=1in]{geometry}

\def\0{\ensuremath\varnothing}

\makeatletter
\newcommand{\compemb}{\mathrel{\mathpalette\comp@emb\relax}}
\newcommand{\comp@emb}[2]{%
  \vcenter{%
    \offinterlineskip\m@th
    \ialign{$#1##$\cr\hookrightarrow\cr\noalign{\vskip1pt}\hookrightarrow\cr}%
  }%
}
\makeatother
\def\a{\ensuremath\alpha }

\def\r{\ensuremath\rho }

\def\g{\ensuremath\gamma }

\def\d{\ensuremath\delta }

\def\f{\ensuremath\varphi }

\def\o{\ensuremath\omega }

\def\w{\ensuremath\bm{w} }

\def\C{\ensuremath\mathcal{C} }

\def\L{\ensuremath\mathcal{L} }

\def\H{\ensuremath\mathcal{H} }

\def\R{\ensuremath\mathbb{R} }

\renewcommand{\div}{\text{div}}

\newlength{\jeroenlen}

\numberwithin{equation}{section}

\title{On the Propulsion of a Rigid Body  in a Viscous Liquid \\ Under the Action of a  Time-Periodic Force} 
 
\author[1]{Giovanni P. Galdi}
\author[2]{Mher M. Karakouzian}

\affil[1]{\small{Correspondence: \href{mailto:galdi@pitt.edu}{galdi@pitt.edu}, Department of Mechanical \& Material Science, Swanson School of Engineering, University of Pittsburgh. Pittsburgh, USA}}
\affil[2]{\small{Department of Mechanical \& Material Science, Swanson School of Engineering, University of Pittsburgh. Pittsburgh, USA}}

\begin{document}
\renewcommand{\labelenumi}{(\roman{enumi})}

    \maketitle

\begin{abstract}
    A rigid body $\mathcal{B}$ moves in an otherwise quiescent viscous liquid filling the whole space outside $\mathcal B$, under the action of a time-periodic force $\bm{\mathsf{f}}$ of period $T$ applied to a given point of $\mathcal B$ and of fixed direction. We assume that the average of $\bm{\mathsf{f}}$ over an interval of length $T$ does not not vanish, and that the amplitude, $\delta$, of $\bm{\mathsf{f}}$ is sufficiently small.  Our goal is to investigate when $\mathcal{B}$ executes a non-zero net motion; that is, $\mathcal{B}$ is able to cover any prescribed distance in a finite time. We show that, at the order $\delta$,  this happens if and only if $\bm{\mathsf{f}}$ and $\mathcal B$ satisfy a certain condition. We also show that this is always the case if $\mathcal B$ is prevented from spinning. Finally, we provide explicit examples where the condition above is satisfied or not. All our analysis is performed in a general class of weak solutions to the coupled system body-liquid problem. \\\\ Keywords: Navier-Stokes, fluid-structure interaction, vibration-induced motion, time-periodic solutions
\end{abstract}

\section{Introduction}\label{ssec:intro}
In recent years, the study of rigid bodies propelling within viscous liquids by means of an applied periodic force has been an active area of research. From the practical point-of-view, this mode of propulsion is advantageous for many reasons. At the macro scale, say for underwater robotics \cite{wilson2006collision, li2023underwater}, it is preferred over the use of fins or propeller blades due to their detrimental effects on the surrounding living organisms and at the micro scale, it provides a primitive means of motion and maneuverability in spatially restrictive environments, such as those found in the human body \cite{hosseini2009design, li2006motion,menciassi2009wireless}, where other means of propulsion are either impractical or impossible.

Most of the research on this subject has been performed either numerically or experimentally by using specific force-producing driving mechanisms (such as moving internal masses or rotors); see \cite{karavaev2016experimental,borisov2018self,vetchanin2013self,borisov2019motion,borisov2020motion} and the reference therein. Concerning a rigorous mathematical study, there are only a few contributions, mainly devoted to  well-posedness of the initial-boundary value problem and the large-time behaviour of the coupled liquid-body system (see \cite{galdi2009motion, galdi2023large, galdi2002strong}). However, the fundamental question of when propulsion actually occurs (specifically, what are the necessary and sufficient conditions) has, to the best of our knowledge, yet to be addressed from a strict mathematical viewpoint. The objective of the current note is to provide a first contribution in this direction.

Specifically, with respect to an inertial frame  $\widehat{\mathcal F}$, consider  a rigid body $\mathcal{B}$ moving in a Navier-Stokes liquid $\L$ that occupies the whole space outside $\mathcal B$. Suppose that at a given point $P$ of the body it is applied  the prescribed force
$$
    \bm{\mathsf{f}}(t)=\mathsf{f}(t)\widehat{\bm{b}},
$$
where $\mathsf{f}$ is a time-periodic function of period $T>0$ (``$T$-periodic") and magnitude $\delta$, and $\widehat{\bm{b}}$ is a constant unit vector. Assume, further, that the force has nonzero average, namely $\overline{\mathsf{f}}\neq 0$ (with the bar denoting the average in time). Our main goal is to find  conditions ensuring that ${\bm{\mathsf{f}}}$ propels $\mathcal{B}$, namely, the center of mass, $G$, of $\mathcal B$ can cover any given distance in a finite time. 
\par
In order to investigate this question, we formulate, as is customary, the governing equations in a  frame attached to $\mathcal{B}$, where  the  domain, $\Omega$, occupied by the liquid becomes time-independent. Thus, denoting by $\mathcal{F}$ such a frame with the origin at $G$, the equations read \cite[Section 1]{galdi2002motion}
\begin{equation}\label{eq:eom_coupled_system}
\begin{aligned}
    \left.\begin{array}{c}
        \displaystyle
        \frac{\partial\bm{v}}{\partial t}+(\bm{v}-\bm{U})\cdot\nabla\bm{v}+\bm{\o}\times\bm{v}=\div\,\textbf{T}(\bm{v},p) \\
        \;\;\;\;\;\;\;\;\;\;\;\;\;\;\;\;\;\;\;\;\;\div\,\bm{v}=0 \\
    \end{array}\right\}&\;\;\;\;\;\;\;{\text{in}\;\Omega\times\mathbb{R}} \\
    \bm{v}=\bm{U}&\;\;\;\;\;\;\;\text{on}\;\partial\Omega\times\mathbb{R} \\
    \lim_{|\bm{x}|\rightarrow\infty}\bm{v}(\bm{x},t)=\textbf{0}\;&\;\;\;\;\;\;\;\text{in}\;\mathbb{R}. \\
        \left.\begin{array}{c}
        \displaystyle
        \;\;M(\dot{\bm{\g}}+\bm{\o}\times\bm{\g})=f\bm{b}-\int_{\partial\Omega}\textbf{T}(\bm{v},p)\cdot\bm{n}\,\text{d}S \\
        \;\;\;\;\;\;\;\;\;\;\;\;\;\;\;\displaystyle\textbf{I}\cdot\dot{\bm{\o}}+\bm{\o}\times(\textbf{I}\cdot\bm{\o})=f(\bm{r}\times\bm{b})-\int_{\partial\Omega}\bm{x}\times(\textbf{T}(\bm{v},p)\cdot\bm{n})\,\text{d}S \\
        \,\displaystyle \frac{d\bm{b}}{dt}=\bm{\o}\times\bm{b} \\
    \end{array}\right\}&\;\;\;\;\;\;\;{\text{in}\;\mathbb{R}}\,.
\end{aligned}
\end{equation}
Here, $\bm{v}$ and $\rho p$ are the velocity and pressure fields of $\mathcal L$, respectively, with $\rho$ its density, whereas $\textbf{T}(\bm{v},p):=-p\textbf{1}+2\nu\textbf{D}(\bm{v})$, with \textbf{1} identity tensor,  $\nu$ kinematic viscosity and $\textbf{D}(\bm{v}):=\frac{1}{2}\left(\nabla\bm{v}+\left(\nabla\bm{v}\right)^{\top}\right)$, is the {Cauchy stress tensor}.  Furthermore, $\rho \,M$  and $\textbf{I}$ represent the mass and the inertia tensor with respect to $G$ of $\mathcal B$, and $\bm{U}:=\bm{\gamma}+\bm{\omega}\times \bm{x}$, with $\bm{\gamma}$ and $\bm{\omega}$  translational and angular velocities of $\mathcal B$, respectively. We also set $f:=\frac{\mathsf{f}}{\r}$ and $\bm{r}:=\overrightarrow{GP}$, while $\bm{n}$ stands for the outer unit normal to $\partial\Omega$.
Finally, the vector $\bm{b}$ is the transformed vector $\widehat{\bm{b}}$ in the body-fixed frame $\mathcal{F}$. More precisely, denoting by $\textbf{Q}=\textbf{Q}(t)$, $t\in \mathbb R$, the one-parameter family of orthogonal matrices associated with the change of frame $\widehat{\mathcal{F}}\rightarrow\mathcal{F}$, we have
\begin{equation}\label{eq:generic_corf}
    \bm{b}(t)=\textbf{Q}^{\top}(t)\cdot\widehat{\bm{b}},\;\;\;\;\;\;\;\textbf{Q}(0):=\textbf{1}.
\end{equation}
Note then that, since the motion of $\mathcal{B}$ is unknown,  so is $\textbf{Q}(t)$ and, therefore, $\bm{b}(t)$. This explains the need for equation $(\ref{eq:eom_coupled_system})_7$ which follows immediately from differentiating (\ref{eq:generic_corf}$)_1$ and then using the property $(\dot{\textbf{Q}}^{\top}\cdot{\textbf{Q}})\cdot\bm{a}=\bm{\o}\times\bm{a}$, $\bm{a}\in\mathbb R^3$.
\par
To address our propulsion problem, we proceed as follows. Since the data, that is, $\sf f$, is $T$-periodic, we look for $T$-periodic weak solutions $(\bm{v},\bm{\gamma},\bm{b})$ to (\ref{eq:eom_coupled_system}). This step is achieved thanks to the results established in \cite{galdi2009motion}, without imposing any restriction on the magnitude $\delta$. Successively, we notice that, denoting by ${\bm{\mathsf{s}}}={\bm{\mathsf{s}}}(t)$ the position of $G$ referred to the frame $\mathcal F$, we have     
\begin{equation}\label{period}
    {\bm{\mathsf{s}}}(T+t)={\bm{\mathsf{s}}}(t)+\int_t^{t+T}\bm{\gamma}(s)\textrm{d}s\,,\ \ \mbox{all $t\in\mathbb R$}\,.
\end{equation}
Therefore, $G$ can cover an arbitrarily given  distance $D$ in a time-span $\tau$ if and only if the average $\overline{\bm{\gamma}}$ of $\bm{\gamma}$ over an interval of length $T$ is not 0. In such a case, from (\ref{period}) it follows that it is enough to take $\tau=NT$, $N\in\mathbb N$,  with $N\ge D/(T|\overline{\bm{\gamma}}|)$. Clearly, $G$ will cover the same distance with respect to the inertial frame $\widehat{\mathcal F}$. As a result, propulsion is reduced to finding conditions on $\bm{f}:=f\,\bm{b}$ and $\mathcal B$ guaranteeing 
\begin{equation}\label{gamma}
    \overline{\bm{\gamma}}\neq\bm{0}.
\end{equation}
With this in mind, we then show that, at first order in $\delta$, equation (\ref{gamma}) is satisfied if and only if
\begin{equation}\label{selfpro}
\widehat{\bm{b}}\neq\textbf{B}\cdot(\bm{r}\times\widehat{\bm{b}})\,,
\end{equation} 
where $\textbf{B}$ is a constant tensor  depending only on the shape of $\mathcal{B}$. Precisely, writing $f=\delta\, F$, 
we prove  that there exists $\delta_0>0$ such that
\begin{equation}\label{selfpro1}
\overline{\bm{\gamma}}=\delta\overline{F}\textbf{A}\cdot\left(\widehat{\bm{b}}-\textbf{B}\cdot(\bm{r}\times\widehat{\bm{b}})\right)+{\sf R}(\delta)\,\ \ \mbox{all $\delta\in (0,\delta_0)$}\,,
\end{equation}
with $\textbf{A}$ a positive definite symmetric tensor and ${\sf R}(\delta)=o(\delta)$ as $\delta\to 0$; see Theorem \ref{thm:suffCond4prop}. This result is obtained by combining the uniform estimates on weak solutions of \cite{galdi2009motion} with a suitable scaling argument in $\delta$. We thus show that, as $\delta\to 0$, the scaled and averaged weak solution must converge, in appropriate topology, to the unique solution of a (time-independent) Stokes problem, for which the associated translational velocity is proved to be non-zero if and only if (\ref{selfpro}) holds. 
\par
It is of some interest to provide examples where (\ref{selfpro1}) is  fulfilled or not. This is done in  Section \ref{ssec:example_sphere} in the simple case where $\mathcal B$ is a sphere. There, we show that (\ref{selfpro1}) holds for any $\widehat{\bm{b}}$ if $\mathcal B$ is homogeneous; otherwise (\ref{selfpro1}) may be violated by choosing $P$ and the location of $G$ appropriately. In this context, it is particularly relevant the situation when $\mathcal B$ is prevented from rotating (for example, by applying a suitable torque on it). In this instance, we formally have $\textbf{B}=\bm{0}$ {\it whatever} the shape of $\mathcal B$ and so, as a result, any $T$-periodic applied force with non-zero average will  propel the body; see Remark~\ref{rem}.    
\par 
The paper is organized as follows. In Section 2, we introduce the function spaces relevant to our problem along with some key estimates. In Section 3, we  give the definition of what it means to be a weak solution to problem (\ref{eq:eom_coupled_system}) and state the corresponding existence result proved in \cite{galdi2009motion}. In the following Section 4 we provide necessary and sufficient conditions for propulsion at the order $\delta$, via the scaling argument mentioned above. Finally, Section 5 is dedicated to the investigation of the validity of (\ref{selfpro}) in the special case where $\mathcal B$ is a sphere.

\hfill

\section{Function Spaces and Related Properties}\label{ssec:fnSpcesand prelim}
We begin to recall some basic notation.
By $B_r$ we indicate the ball of radius $r>0$ in $\mathbb{R}^3$ centered at the origin and set $\mathbb{S}^2:=\partial B_1$. For a domain $A\subseteq\mathbb{R}^3$,   $L^q(A)$  denotes the usual Lebesgue space endowed with the norm $\|\cdot\|_{L^q(A)}$ and, for $m\in\mathbb{N}$ and $q\in[1,\infty]$, $W^{m,q}(A)$ stands for the Sobolev space with norm $\|\cdot\|_{W^{m,q}(A)}$. Moreover, $D^{m,q}(A)$ will denote the homogeneous Sobolev space with semi-norm $|u|_{D^{m,q}(A)}:=\sum_{m=|\a|}\|D^{\a}u\|_{L^q(A)}$. When any of the above function spaces are used with the subscript ``per", we shall mean that a function $u$ of this space has the additional property of being $T$-periodic; namely, $u(t+T)=u(t)$, for all $t\in\mathbb{R}$. Finally, for a function $w=w(t)$ defined in the interval $(0,T)$ we define the average:
$$
    \overline{w}:=\frac1T\int_0^Tw(t)\,\textrm{d}t.
$$
\par
Let $\Omega\subset\mathbb{R}^3$ be a locally Lipschitz exterior domain, where $\Omega:=\mathbb{R}^3\setminus\overline{\mathcal{B}}$, for some bounded domain $\mathcal{B}\subset \mathbb{R}^3$. Physically, $\mathcal{B}$ is the moving rigid body described in Section \ref{ssec:intro}. For each $R>\text{diam}\,\mathcal{B}$, we also adopt the following convention:
$$
    \Omega_R:=\Omega\cap B_R.
$$

Now, define
$$
    \mathcal{R}:=\{\bm{U}\in C^{\infty}(\mathbb{R}^3):\bm{U}(\bm{x})=\bm{\g}+\bm{\o}\times\bm{x},\;\text{for some $\bm{\g},\bm{\o}\in \mathbb{R}^3$}\}
$$
and if $\bm{U}\in\mathcal{R}$ depends on a vector $\bm{v}$, let us write $\bm{U}_{\bm{v}}$ and define the vectors $\bm{\g}_{\bm{v}}, \bm{\o}_{\bm{v}}\in\mathbb{R}^3$ to be those which correspond to $\bm{U}_{\bm{v}}$ in the definition of $\mathcal{R}$; that is,
$$
    \bm{U}_{\bm{v}}=\bm{\g}_{\bm{v}}+\bm{\o}_{\bm{v}}\times\bm{x}.
$$
For $A\in \{\Omega,\Omega_R\}$ we introduce the set
\begin{equation*}
    \begin{aligned}
    \C(A)&:= \left\{\bm{\f}\in C_0^{\infty}(\overline{A}):
    \begin{array}{l}
    \text{$\div\,\bm{\f}=0$ in $A$;} \\
    \text{$\bm{\f}=\bm{U}_{\bm{\f}}$\  in a neighborhood of $\overline{\mathcal B}$, for some $\bm{U}_{\bm{\f}} \in \mathcal{R}$;} \\
    \text{$\bm{\f}=\bm{0}$ in a neighborhood of $\partial B_R$ if $A\equiv \Omega_R$}
    \end{array}
    \right\},
    \end{aligned}
\end{equation*}
and define the inner product:
$$
    \left(\bm{u},\w\right)_{\H(A)}:=\int_{A}\textbf{D}(\bm{u}):\textbf{D}(\w)\,\text{d}V,\;\;\;\text{for all}\;\bm{u},\w\in \C(A),
$$
with associated norm:
$$
    \|\bm{u}\|_{\H(A)}:=\|\textbf{D}(\bm{u})\|_{L^2(A)},\;\;\;\text{for all}\;\bm{u}\in \C(A),
$$
respectively. Finally,  we set
$$
\H(A):=\overline{\C(A)}^{\|\cdot\|_{\H(A)}}.
$$
It can be shown (see   \cite[Lemma 4.11]{galdi2002motion}), that
\begin{equation*}
    \begin{aligned}
    \H(\Omega)&:= \left\{\bm{v}\in W^{1,2}_\text{loc}(\mathbb R^3)\cap L^6(\mathbb R^3):
    \begin{array}{l}
    \text{$\textbf{D}(\bm{v})\in L^2(\mathbb R^3)$, $\div\,\bm{v}=0$, and} \\
    \text{$\bm{v}=\bm{U}_{\bm{v}}$ in $\overline{\mathcal B}$, for some $\bm{U}_{\bm{v}} \in \mathcal{R}$ }
    \end{array}
    \right\}.
    \end{aligned}
\end{equation*}
Likewise, for the ``local space", the following characterization holds:
\begin{equation*}
    \H(\Omega_R)  = \{\bm{v}\in W^{1,2}(B_R): \text{$\div\,\bm{v}=0$ in $\Omega_R$; $\bm{v}=\bm{U}_{\bm{v}}$ in $\overline{\mathcal B}$; $\bm{v}=\bm{0}$ around $\partial B_R$}\}.
\end{equation*}
It is known that $\H(A)$ is a Hilbert space with the norm $(\cdot,\cdot)_{\H(A)}$. For $m\in\mathbb{N}\cup\{\infty\}$ and fixed period $T>0$, we introduce the test function spaces
\begin{equation*}
    \begin{aligned}
    \C^m_\text{per}(A\times\mathbb{R})&:= \left\{\bm{\f}\in C^m(A\times\mathbb{R}):
    \begin{array}{l}
    \text{$\div\,\bm{\f}=0$ in $A$; $\bm{\f}$ is $T$-periodic;} \\
    \text{there exists $\bm{U}_{\bm{\f}} \in C_\text{per}^m(\mathbb{R};\mathcal{R})$, such that} \\
    \text{$\bm{\f}(\bm{x},\cdot)=\bm{U}_{\bm{\f}}(\bm{x},\cdot)$, for all $\bm{x}\in\overline{\mathcal B}$;} \\
    \text{there exists $r>\text{diam}\,\mathcal{B}$, such that $\bm{\f}(\bm{x},t)=0$,} \\
    \text{for all $\bm{x}\in\mathbb{R}^3\setminus\overline{B}_{r}$ and all $t\in\mathbb{R}$, where $r<R$ if $A\equiv \Omega_R$}
    \end{array}
    \right\},
    \end{aligned}
\end{equation*}
where we use $\C^m_\text{per}(A\times[0,T])$ to denote the functions of $\C^m_\text{per}(A\times\mathbb{R})$ restricted to $[0,T]$. Similarly, we will use $C^m_\text{per}([0,T])$ to denote the functions of $C^m_\text{per}(\mathbb{R})$ restricted to $[0,T]$.

We conclude this section with the following lemma, containing a collection of important estimates pertaining to the space $\H(A)$ (see  \cite[Section 4]{galdi2002motion}).

\hfill

\begin{lemma}\label{thm:KornsIdentity_coupledSyst}
    For $R>\text{diam}\,\mathcal{B}$, let $A\in\{\Omega,\Omega_R\}$ and $\bm{u}\in\H(A)$. Then
    \begin{equation}\label{eq:KornsIdentity_coupledSyst}
        \|\nabla \bm{u}\|_{L^2(A)}=\sqrt{2}\|\bm{u}\|_{\H(A)}
    \end{equation}
    and there exist $c_1,c_2>0$, independent of $A$, such that, for all $\bm{u}\in\H(A)$, the following inequalities hold:
    \begin{eqnarray}
        |\bm{\g}_{\bm{u}}|+|\bm{\o}_{\bm{u}}|&\leq& c_1\|\bm{u}\|_{\H(A)};\label{eq:estOnGamma_u} \\
        \|\bm{u}\|_{L^6(A)} &\leq& c_2 \|\bm{u}\|_{\H(A)}.\label{eq:sobolevInHO_R}
    \end{eqnarray}
\end{lemma}

\hfill

\section{\textbf{Weak Solutions to the Coupled Liquid-Body Problem}}\label{WSBL}
Let us begin by furnishing a weak formulation for problem (\ref{eq:eom_coupled_system}). Formally dot-multiplying (\ref{eq:eom_coupled_system}$)_1$ by arbitrary $\bm{\f}\in \C_\text{per}^1(\Omega\times\mathbb{R})$ and integrating by parts using (\ref{eq:eom_coupled_system}$)_\text{2,3,5,6}$ and also periodicity, we get
\begin{equation}\label{eq:wkForm_coupledSyst-1}
    \begin{aligned}
            &\int_0^T {\Bigg[}\left(\bm{v},\frac{\partial\bm{\f}}{\partial t}\right)_{L^2(\Omega)}+M\bm{\xi}\cdot\dot{\bm{\g}}_{\bm{\f}}+\bm{\o}\cdot\textbf{I}\cdot\dot{\bm{\o}}_{\bm{\f}}-\left((\bm{v}-\bm{U})\cdot\nabla\bm{v},\bm{\f}\right)_{L^2(\Omega)}-2\nu\left( \bm{v},\bm{\f}\right)_{\H(\Omega)}\\
            &\;\;\;\;\;\;\;\;\;\;\;\;\;\;\;\;\;\;\;\;\;\;\;\;\;\;\;-M(\bm{\o}\times\bm{\g})\cdot\bm{\g}_{\bm{\f}}-[\bm{\o}\times(\textbf{I}\cdot\bm{\o})]\cdot\bm{\o}_{\bm{\f}}+f\bm{b}\cdot\bm{\g}_{\bm{\f}}+f(\bm{r}\times\bm{b})\cdot\bm{\o}_{\bm{\f}}{\Bigg]}\text{d}t=0.
    \end{aligned}
    \end{equation}
Similarly, multiplying (\ref{eq:eom_coupled_system}$)_7$ by arbitrary $\bm{\psi}\in C^1_\text{per}([0,T])$ and integrating by parts, we get
\begin{equation}\label{eq:wkForm_coupledSyst-2}
    \int_0^T\left[\bm{b}\cdot\dot{\bm{\psi}}+(\bm{\o}\times\bm{b})\cdot\bm{\psi}\right]\,\text{d}t=0.
\end{equation}
Then, as in \cite{galdi2009motion}, we give the following definition.

\hfill

\begin{definition}\label{def:wkSol-1}
    Let $f\in L_\text{per}^{\infty}(\mathbb{R})$. Then, $(\bm{v}, \bm{\g},\bm{\o},\bm{b})$ is said to be a {\em $T$-periodic weak solution} to problem (\ref{eq:eom_coupled_system}) if
\begin{enumerate}
    \item $\bm{v}\in L_\text{per}^2(\mathbb{R};\H(\Omega))$ and $\bm{\g},\bm{\o}\in L^2_\text{per}(\mathbb{R})$ with $\bm{v}=\bm{U}:=\bm{\g}+\bm{\o}\times\bm{x}$ in $\overline{\mathcal B}$\,\footnote{In fact, due to (\ref{eq:estOnGamma_u}), the condition $\bm{\g},\bm{\o}\in L^2_\text{per}(\mathbb{R})$ is automatic if we take $\bm{U}\equiv\bm{U}_{\bm{v}}$.};
    \item $\bm{b}\in W^{1,2}_\text{per}(\mathbb{R};\mathbb{S}^2)$
    \item $(\bm{v}, \bm{\g},\bm{\o},\bm{b})$ verifies (\ref{eq:wkForm_coupledSyst-1}) for all $\bm{\f}\in\C^1_\text{per}(\Omega\times \mathbb{R})$ and (\ref{eq:wkForm_coupledSyst-2}) for all $\bm{\psi}\in C^1_\text{per}([0,T])$.
\end{enumerate}
\end{definition}

\hfill

In the sense of the above definition, existence of weak solutions to problem (\ref{eq:eom_coupled_system}) has been shown in \cite[Theorem 3.4]{galdi2009motion} along with appropriate estimates. These results are summarized in the following theorem.

\hfill

\begin{thm}\label{thm:exist_wk_sol_coupledSyst}
    Let $f\in L_\text{per}^{\infty}(\mathbb{R})$. There exists a $T$-periodic weak solution $(\bm{v},\bm{\g},\bm{\o},\bm{b})$ to problem (\ref{eq:eom_coupled_system}). Moreover, there is a constant $C=C(T,\nu,\bm{r})>0$, such that
    \begin{equation}\label{eq:energyForWRFinal_coupledSyst}
        \|\bm{v}\|_{L^2(0,T;\H(\Omega))}+\|\bm{\gamma}\|_{L^2(0,T)}+\|\bm{\omega}\|_{L^2(0,T)}+\left|\left|\frac{d\bm{b}}{dt}\right|\right|_{L^2(0,T)}\leq C\|f\|_{L^{\infty}(0,T)}.
    \end{equation}
\end{thm}

\hfill

\section{\textbf{Sufficient Conditions for Propulsion}}\label{SCFP}
For each $\d>0$, consider the following scaled decompositions of the vector fields $\bm{v}$, $\bm{\g}$, and $\bm{\o}$ from Theorem \ref{thm:exist_wk_sol_coupledSyst}:
\begin{equation}\label{eq:decomp_1stappearance}
    \bm{v}=\d(\bm{u}+\w),\;\;\;\;\;\;\;\;\;\;\bm{\g}=\d(\bm{\xi}+\bm{\chi}),\;\;\;\;\;\text{and}\;\;\;\;\;\bm{\o}=\d(\bm{\zeta}+\bm{\eta}),
\end{equation}
where $\delta\bm{u}:=\overline{\bm{v}}$, $\delta\bm{\xi}:=\overline{\bm{\g}}$, $\delta\bm{\zeta}:=\overline{\bm{\o}}$, $\delta\w:=\bm{v}-\overline{\bm{v}}$, $\delta\bm{\chi}:=\bm{\g}-\overline{\bm{\g}}$, and $\delta\bm{\eta}:=\bm{\o}-\overline{\bm{\o}}$. The vectors $\bm{u}$, $\bm{\xi}$, and $\bm{\zeta}$ are then the (scaled) \textit{time-averaged components} of $\bm{v}$, $\bm{\g}$, and $\bm{\o}$, respectively, with $\w$, $\bm{\chi}$, and $\bm{\eta}$ their respective (scaled) \textit{purely oscillatory components}. Consequently, these components satisfy
\begin{equation}\label{eq:props_of_components}
    \frac{\partial\bm{u}}{\partial t}=\frac{\partial\bm{\xi}}{\partial t}=\frac{\partial\bm{\zeta}}{\partial t}=\overline{\w}=\overline{\bm{\chi}}=\overline{\bm{\eta}}=\textbf{0}.
\end{equation}
Let us also scale the force, say $\d F:=f$. Then, substituting these expressions into (\ref{eq:eom_coupled_system}$)_\text{1-3,5-6}$, taking the average over $(0,T)$, and using the properties of $\bm{u}$, $\bm{\xi}$, $\w$, and $\bm{\chi}$ above, we get
\begin{equation}\label{eq:eom_avged-nonLnr-coupled_system}
\begin{aligned}
    \left.\begin{array}{c}
        \displaystyle
        \delta\left[(\bm{u}-\bm{\xi}-\bm{\zeta}\times\bm{x})\cdot\nabla\bm{u}+\overline{(\w-\bm{\chi}-\bm{\eta}\times\bm{x})\cdot\nabla\w}\right]=\div\textbf{T}(\bm{u},\pi) \\
        \;\;\;\;\;\;\;\;\;\;\;\;\;\;\;\;\;\;\;\;\;\;\;\;\;\;\;\;\;\;\;\;\;\;\;\;\;\;\;\;\;\;\;\;\;\;\;\;\;\;\;\;\;\;\;\div\,\bm{u}=0 \\
    \end{array}\right\}&\;\;\;\;\;\;\;{\text{in}\;\Omega} \\
    \bm{u}=\bm{\xi}+\bm{\zeta}\times\bm{x}&\;\;\;\;\;\;\;\text{on}\;\partial\Omega \\
    \d M\left[\bm{\zeta}\times\bm{\xi}+\overline{\bm{\eta}\times\bm{\chi}}\right]=\overline{F\bm{b}}-\int_{\partial\Omega} \textbf{T}(\bm{u},\pi)\cdot\bm{n}\;\text{d}S,& \\
    \d \left[\bm{\zeta}\times(\textbf{I}\cdot\bm{\zeta})+\overline{\bm{\eta}\times(\textbf{I}\cdot\bm{\eta})}\right]=\bm{r}\times\overline{F\bm{b}}-\int_{\partial\Omega} \bm{x}\times\textbf{T}(\bm{u},\pi)\cdot\bm{n}\;\text{d}S,&
\end{aligned}
\end{equation}
where $\delta\pi:=\overline{p}$, and substituting (\ref{eq:decomp_1stappearance}$)_3$ into (\ref{eq:eom_coupled_system}$)_7$, we get
\begin{equation}\label{eq:separate_eqn_for-b}
    \frac{d\bm{b}}{dt}=\d(\bm{\zeta}+\bm{\eta})\times\bm{b}.
\end{equation}
Formally taking $\delta\rightarrow 0$ in equation (\ref{eq:separate_eqn_for-b}), we see that $\bm{b}$ tends to some constant vector $\bm{b}_0\in\mathbb{R}^3$. In fact, from (\ref{eq:generic_corf}), apparently $\bm{b}_0=\widehat{\bm{b}}$ (we shall soon make this precise). Then, in the limit $\delta\rightarrow 0$, from (\ref{eq:eom_avged-nonLnr-coupled_system}) we (formally) obtain the following (time-independent) Stokes problem:
\begin{equation}\label{eq:eom_avged-Lnr-coupled_system}
\begin{aligned}
    \left.\begin{array}{c}
        \displaystyle \div\textbf{T}(\bm{u}_0,\pi_0)=\bm{0} \\
        \;\,\;\;\;\;\;\;\;\;\div\,\bm{u}_0=0 \\
    \end{array}\right\}&\;\;\;\;\;\;\;{\text{in}\;\Omega} \\
    \bm{u}_0=\bm{\xi}_0+\bm{\zeta}_0\times\bm{x}&\;\;\;\;\;\;\;\text{on}\;\partial\Omega \\
    \int_{\partial\Omega} \textbf{T}(\bm{u}_0,\pi_0)\cdot\bm{n}\;\text{d}S=\overline{F}\widehat{\bm{b}},& \\
    \int_{\partial\Omega} \bm{x}\times\textbf{T}(\bm{u}_0,\pi_0)\cdot\bm{n}\;\text{d}S=\bm{r}\times\overline{F}\widehat{\bm{b}}.&
\end{aligned}
\end{equation}
Again formally multiplying (\ref{eq:eom_avged-Lnr-coupled_system}$)_1$ by arbitrary $\bm{\psi}\in\H(\Omega)$ and integrating by parts over $\Omega$, as was done to obtain (\ref{eq:wkForm_coupledSyst-1}), we are lead to a weak formulation of (\ref{eq:eom_avged-Lnr-coupled_system}), made precise by the following definition.

\hfill

\begin{definition}\label{def:wkSol-2}
    Let $F\in L_\text{per}^{\infty}(\mathbb{R})$. Then $(\bm{u}_0,\bm{\xi}_0,\bm{\zeta}_0)$ is a {\em weak solution} to the Stokes problem (\ref{eq:eom_avged-Lnr-coupled_system}) if
    \begin{enumerate}
        \item $\bm{u}_0\in \H(\Omega)$ and $\bm{\xi}_0,\bm{\zeta}_0\in\mathbb{R}^3$ are such that $\bm{u}_0=\bm{\xi}_0+\bm{\zeta}_0\times\bm{x}$ on $\partial\Omega$;
        \item $\bm{u}_0$ satisfies
        \begin{equation}\label{eq:wkForm_lnrProb}
            2\nu\left(\bm{u}_0,\bm{\psi}\right)_{\H(\Omega)}=\overline{F}\widehat{\bm{b}}\cdot\bm{\g}_{\bm{\psi}}+\bm{r}\times\overline{F}\widehat{\bm{b}}\cdot\bm{\o}_{\bm{\psi}},\;\;\;\;\;\text{for every $\bm{\psi}\in \H(\Omega)$}.
        \end{equation}
    \end{enumerate}
\end{definition}

\hfill

Now, for each $\d>0$, thanks to Theorem \ref{thm:exist_wk_sol_coupledSyst}, we have a weak solution $(\bm{v}_{\delta},\bm{\g}_{\delta},\bm{\o}_{\delta},\bm{b}_{\delta})$ to problem (\ref{eq:eom_coupled_system}). We claim that, as $\delta\rightarrow 0$, the vector fields $\bm{b}_{\delta}$ converge (in some suitable sense) to $\widehat{\bm{b}}$, meanwhile the corresponding time-averaged parts $(\bm{u}_{\delta},\bm{\xi}_{\delta},\bm{\zeta}_{\delta})$ converge to the weak solutions $(\bm{u}_0,\bm{\xi}_0,\bm{\zeta}_0)$ of (\ref{eq:eom_avged-Lnr-coupled_system}), whose both properties of existence and uniqueness must  be verified first. To this end, we recall the following result, for whose proof we refer to \cite[Section V.4]{galdi2011introduction},   \cite[Sections 5.2-5.4]{happel1983low}. \\

\begin{lemma}\label{thm:eom_aux}
    Let $s\in(1,\infty)$, $q\in \left(\frac{3}{2},\infty\right)$ and $r\in (3,\infty)$. For each $i=1,2,3$, there exists unique solutions 
$$(\bm{h}^{(i)},p^{(i)}),(\bm{H}^{(i)},P^{(i)})\in [D^{2,s}(\Omega)\cap D^{1,q}(\Omega)\cap L^r(\Omega)\cap C^\infty(\Omega)]\times [D^{1,s}(\Omega)\cap L^q(\Omega)\cap C^\infty(\Omega)]
$$ 
to the Stokes problems
    \begin{equation}\label{eq:eom_aux-1}
    \begin{aligned}
    \left.\begin{array}{c}
        \displaystyle
        \div\,\textbf{T}(\bm{h}^{(i)},p^{(i)})=\bm{0} \\
        \;\;\;\;\;\;\;\;\;\;\;\,\div\,\bm{h}^{(i)}=0 \\
    \end{array}\right\}&\;\;\;\;\;\;\;{\text{in}\;\Omega} \\
    \bm{h}^{(i)}=\bm{e}_i&\;\;\;\;\;\;\;\text{on}\;\partial\Omega
\end{aligned}
\end{equation}
and
    \begin{equation}\label{eq:eom_aux-2}
    \begin{aligned}
    \left.\begin{array}{c}
        \displaystyle
        \div\,\textbf{T}(\bm{H}^{(i)},P^{(i)})=\bm{0} \\
        \;\;\;\;\;\;\;\;\;\;\;\;\,\div\,\bm{H}^{(i)}=0 \\
    \end{array}\right\}&\;\;\;\;\;\;\;{\text{in}\;\Omega} \\
    \bm{H}^{(i)}=\bm{e}_i\times\bm{x}&\;\;\;\;\;\;\;\text{on}\;\partial\Omega.
\end{aligned}
\end{equation}
Moreover, for $i,k\in\{1,2,3\}$, defining (component-wise) the matrices
\begin{equation}
\begin{aligned}\label{matrices}
    (\textbf{K})_{ki}:=\bm{e}_k\cdot\int_{\partial\Omega}\left(\textbf{T}(\bm{h}^{(i)},p^{(i)})\cdot\bm{n}\right) \text{d}S,  &\;\;\;\;\;\;\;\;\;\;(\textbf{C})_{ki}:=\bm{e}_k\cdot\int_{\partial\Omega}\left(\bm{x}\times\textbf{T}(\bm{h}^{(i)},p^{(i)})\cdot\bm{n}\right) \text{d}S, \\
    (\textbf{S})_{ki}:=\bm{e}_k\cdot\int_{\partial\Omega}\left(\textbf{T}(\bm{H}^{(i)},P^{(i)})\cdot\bm{n}\right) \text{d}S,  &\;\;\;\;\;\;\;\;\;\;(\bm{\Uptheta})_{ki}:=\bm{e}_k\cdot\int_{\partial\Omega}\left(\bm{x}\times\textbf{T}(\bm{H}^{(i)},P^{(i)})\cdot\bm{n}\right) \text{d}S,
\end{aligned}
\end{equation}
we have that $\textbf{K}$ and $\bm{\Uptheta}$ are both symmetric and invertible and $\textbf{S}=\textbf{C}^{\top}$. Finally, both matrices $\textbf{K}-\textbf{C}\cdot\bm{\Uptheta}^{-1}\cdot\textbf{C}^{\top}$ and $\bm{\Uptheta}-\textbf{C}^{\top}\cdot\textbf{K}^{-1}\cdot\textbf{C}$ are invertible as well.
\end{lemma}

\hfill

Observe that, for each $i=1,2,3$, problems (\ref{eq:eom_aux-1}) and (\ref{eq:eom_aux-2}) describe the flow of a viscous liquid around a body with the \textit{prescribed} motion of pure translation along basis vector $\bm{e}_i$ for (\ref{eq:eom_aux-1}) and of pure rotation about the axis directed along $\bm{e}_i$ for (\ref{eq:eom_aux-2}). In turn, $(\textbf{K})_{ki}$ represents the $k^\text{th}$ component of the hydrodynamic force exerted on $\partial\Omega$ due to \textit{pure translation} along the direction $\bm{e}_i$ and each $(\textbf{S})_{ki}$ represents those due to \textit{pure rotation} about the axis directed along $\bm{e}_i$. Analogously, the components of $\textbf{C}$ and $\bm{\Uptheta}$ represent the hydrodynamic torques with respect to $G$, due to pure translation and pure rotation, respectively.

\begin{lemma}\label{thm:eom_steady_lowRe2}
    For any given $F\in L_\text{per}^{\infty}(\mathbb{R})$, there exists a unique corresponding weak solution to problem (\ref{eq:eom_avged-Lnr-coupled_system}) satisfying, in addition,  
$$
(\bm{u}_0,\bm{\xi}_0,\bm{\zeta}_0)\in [D^{2,s}(\Omega)\cap D^{1,q}(\Omega)\cap L^r(\Omega)]\times \mathbb R^3\times\mathbb R^3\,,\ \ s\in(1,\infty), \ q\in(\mbox{$\frac32$},\infty),\ r\in (3,\infty)\,.
$$ 
Furthermore, $\bm{u}_0\in C^\infty(\Omega)$, and there exists $\pi_0\in C^\infty(\Omega)\cap D^{1,s}(\Omega)\cap L^q(\Omega)$ such that $(\bm{u}_0,p_0,\bm{\xi}_0,\bm{\zeta}_0)$ solves (\ref{eq:eom_avged-Lnr-coupled_system}) in the ordinary sense.
\end{lemma}
\begin{proof}
Let
\begin{equation}\label{def:wkSoltoStokesProb-1}
\begin{aligned}
    \bm{\xi}_0&:=\left(\textbf{K}-\textbf{C}\cdot\bm{\Uptheta}^{-1}\cdot\textbf{C}^{\top}\right)^{-1}\cdot\left(\overline{F}\widehat{\bm{b}}-\textbf{C}\cdot\bm{\Uptheta}^{-1}\cdot(\bm{r}\times\overline{F}\widehat{\bm{b}})\right) \\
    \bm{\zeta}_0&:=\left(\bm{\Uptheta}-\textbf{C}^{\top}\cdot\textbf{K}^{-1}\cdot\textbf{C}\right)^{-1}\cdot\left(\bm{r}\times\overline{F}\widehat{\bm{b}}-\textbf{C}^{\top}\cdot\textbf{K}^{-1}\cdot\overline{F}\widehat{\bm{b}}\right)
\end{aligned}
\end{equation}
and, for $\bm{\xi}_0=\xi_{0i}\bm{e}_i$ and $\bm{\zeta}_0=\zeta_{0i}\bm{e}_i$, define
\begin{equation}\label{def:wkSoltoStokesProb-2}
    \bm{u}_0:=\sum_{i=1}^3\left(\xi_{0i}\bm{h}^{(i)}+\zeta_{0i}\bm{H}^{(i)}\right)\;\;\;\;\;\;\;\;\;\;\text{and}\;\;\;\;\;\;\;\;\;\;\pi_0:=\sum_{i=1}^3\left(\xi_{0i}p^{(i)}+\zeta_{0i}P^{(i)}\right).
\end{equation}
In view of  Lemma \ref{thm:eom_aux}, we infer that $(\bm{u}_0,\bm{\xi}_0, \bm{\zeta}_0, \pi_0)$, possesses all the stated regularity properties. Furthermore, multiplying (\ref{eq:eom_aux-1}$)_\text{1-3}$ by $\xi_{0i}$ and summing over $i$,  then multiplying (\ref{eq:eom_aux-2}$)_\text{1-3}$ by $\zeta_{0i}$ and adding the resulting equations, we immediately obtain that $(\bm{u}_0,\bm{\xi}_0, \bm{\zeta}_0, \pi_0)$ satisfies (\ref{eq:eom_avged-Lnr-coupled_system}$)_\text{1-3}$. Next, solving for $\overline{F}\widehat{\bm{b}}$ and $\bm{r}\times\overline{F}\widehat{\bm{b}}$ in (\ref{def:wkSoltoStokesProb-1}), we get
\begin{equation}\label{def:wkSoltoStokesProb-3}
\begin{aligned}
    \overline{F}\widehat{\bm{b}} &=\textbf{K}\cdot\bm{\xi}_0+\textbf{C}\cdot\bm{\zeta}_0 \\
    \bm{r}\times\overline{F}\widehat{\bm{b}} &=\textbf{C}^{\top}\cdot\bm{\xi}_0+\bm{\Uptheta}\cdot\bm{\zeta}_0.
\end{aligned}
\end{equation}
Employing (\ref{def:wkSoltoStokesProb-3}$)_1$ in combination with Lemma \ref{thm:eom_aux}, one easily verifies also the validity of (\ref{eq:eom_avged-Lnr-coupled_system}$)_4$ and, similarly, from (\ref{def:wkSoltoStokesProb-3}$)_2$ one obtains (\ref{eq:eom_avged-Lnr-coupled_system}$)_5$, thus completing the proof of existence.
Concerning uniqueness, let $\bm{u}_0'$ be another weak solution to (\ref{eq:eom_avged-Lnr-coupled_system}) in the sense of (i)-(ii) of the definition. Then, $\bm{u}_0'$ satisfies
    $$
        2\nu\left(\bm{u}_0',\bm{\psi}\right)_{\H(\Omega)}=\overline{F}\widehat{\bm{b}}\cdot\bm{\g}_{\bm{\psi}}+\bm{r}\times\overline{F}\widehat{\bm{b}}\cdot\bm{\o}_{\bm{\psi}},\;\;\;\;\;\text{for every $\bm{\psi}\in \H(\Omega)$}.
    $$
    The result then follows by subtracting this from (\ref{eq:wkForm_lnrProb}) and taking, in particular, $\bm{\psi}:=\bm{u}_0-\bm{u}_0'$.
\end{proof}

We are now in a position to prove the convergences claimed earlier on.

\hfill

\begin{lemma}\label{thm:thrustApprox}
    Let $F\in L_\text{per}^{\infty}(\mathbb{R})$ and $\d>0$. Let $(\bm{v}_{\delta},\bm{\g}_{\delta},\bm{\o}_{\delta},\bm{b}_{\delta})$ be a weak solution to problem (\ref{eq:eom_coupled_system}) corresponding to $f_{\delta}:=\d F$ and apply the decomposition from (\ref{eq:decomp_1stappearance}) to $\bm{v}_{\delta}$ and $\bm{\g}_{\delta}$:
    \begin{equation*}
    \bm{v}_{\delta}=\d(\bm{u}_{\delta}+\w_{\delta})\;\;\;\;\;\;\;\;\;\;\bm{\g}_{\delta}=\d(\bm{\xi}_{\delta}+\bm{\chi}_{\delta})\;\;\;\;\;\;\;\;\;\;\bm{\o}_{\delta}=\d(\bm{\zeta}_{\delta}+\bm{\eta}_{\delta}).
    \end{equation*}
    Then, as $\delta\rightarrow 0$,
    \begin{equation}\label{eq:PropConvgces}
    \begin{aligned}
        \bm{b}_{\delta}&\longrightarrow \widehat{\bm{b}}\;\;\;\;\text{in $C([0,T];\mathbb{S}^2)$} \\
        \bm{u}_{\delta} & \xrightharpoonup{\;\;\;\;\,} \bm{u}_0\;\;\,\text{in $\H(\Omega)$}, \\
        \bm{\xi}_{\delta}&\longrightarrow \bm{\xi}_0\;\;\;\text{in $\mathbb{R}^3$},\;\;\;\;\;\text{and} \\
        \bm{\zeta}_{\delta}&\longrightarrow \bm{\zeta}_0\;\;\;\text{in $\mathbb{R}^3$}, \\
    \end{aligned}
    \end{equation}
    where $(\bm{u}_0,\bm{\xi}_0,\bm{\zeta}_0)$ is the weak solution to problem (\ref{eq:eom_avged-Lnr-coupled_system}) furnished by Lemma \ref{thm:eom_steady_lowRe2}.
\end{lemma}
\begin{proof}
    By the uniqueness property afforded by Lemma \ref{thm:eom_steady_lowRe2}, it suffices to show (\ref{eq:PropConvgces}) for a subsequence $\{\d_n\}_{n\in\mathbb{N}}$, say, of strictly positive numbers with $\lim_{n\rightarrow\infty}\d_n=0$. Given such a sequence, write
    \begin{equation}\label{eq:expr_4_fn}
        f_n:=\d_n F
    \end{equation}
    and, for each $n\in\mathbb{N}$, let $(\bm{v}_n,\bm{\g}_n, \bm{\o}_n,\bm{b}_n)$ be a weak solution to problem (\ref{eq:eom_coupled_system}) corresponding to $f_n$. As in the theorem statement, also write,
    \begin{equation}\label{eq:scaledQties-2}
    \bm{v}_n=\d_n(\bm{u}_n+\w_n)\;\;\;\;\;\;\;\;\;\;\bm{\g}_n=\d_n(\bm{\xi}_n+\bm{\chi}_n)\;\;\;\;\;\;\;\;\;\;\bm{\o}_n=\d_n(\bm{\zeta}_n+\bm{\eta}_n).
    \end{equation}
    First, substituting (\ref{eq:expr_4_fn}) in (\ref{eq:energyForWRFinal_coupledSyst}) and passing to the limit as $n\rightarrow\infty$, we immediately deduce
    \begin{equation}\label{eq:bd_4_bn}
        \lim_{n\rightarrow\infty}\left\|\frac{d\bm{b}_n}{dt}\right\|_{L^2(0,T)}=0.
    \end{equation}
    Since also $\bm{b}_n\in\mathbb{S}^2$, we have that $\bm{b}_n$ is bounded uniformly in $W^{1,2}(0,T)$, and so, by elementary embedding inequality,
\begin{equation}\label{eq:cvgce_of_bntob_0}
    \bm{b}_n\longrightarrow\bm{b}_0\;\;\;\;\text{in $C([0,T];\mathbb{R}^3)$},
\end{equation}  
for some $\bm{b}_0$. Then, we can pass to the limit in the property $|\bm{b}_n|=1$, to conclude that (\ref{eq:cvgce_of_bntob_0}) holds with $\mathbb{S}^2$ in place of $\mathbb{R}^3$. Furthermore, by  (\ref{eq:bd_4_bn}) and dominated convergence, we have that $\bm{b}_0$ is a constant, so upon passing to the limit as $n\to\infty$ in  (\ref{eq:generic_corf}) with $\bm{b}\equiv\bm{b}_n$, and employing (\ref{eq:cvgce_of_bntob_0}), we conclude that $\bm{b}_0\equiv \widehat{\bm{b}}$, thus proving (\ref{eq:PropConvgces}$)_1$.
Next, thanks to  Theorem \ref{thm:exist_wk_sol_coupledSyst} 
and the definition of $(\bm{u}_n,\bm{\xi}_n,\bm{\zeta}_n)$, we obtain the estimate
    \begin{equation}\label{eq:est_on_un}
        \|\bm{u}_n\|_{\H(\Omega)}+|\bm{\xi}_n|+|\bm{\zeta}_n|\leq {\kappa}\,,
    \end{equation}
where, from now on, by $\kappa$ we denote a generic positive constant depending, at most, on $F$ and $T$.    
Then, by standard compactness theorems, one can find $\widetilde{\bm{u}}\in\H(\Omega)$ and $\widetilde{\bm{\xi}},\widetilde{\bm{\zeta}}\in\mathbb{R}^3$, such that (up to subsequence)
    \begin{equation}\label{eq:propConvergences}
        \bm{u}_n \xrightharpoonup{\;\;\;\;\,} \widetilde{\bm{u}}\;\;\;\text{in $\H(\Omega)$},\;\;\;\;\;\bm{\xi}_n\longrightarrow \widetilde{\bm{\xi}}\;\;\;\text{in $\mathbb{R}^3$},\;\;\;\;\;\text{and}\;\;\;\;\;\bm{\zeta}_n\longrightarrow \widetilde{\bm{\zeta}}\;\;\;\text{in $\mathbb{R}^3$} 
    \end{equation}
    as $n\rightarrow\infty$ and such that these limits satisfy, in particular, the boundary condition
    \begin{equation}\label{eq:trace-result-in-prop-thm}
        \widetilde{\bm{u}}=\widetilde{\bm{\xi}}+\widetilde{\bm{\zeta}}\times\bm{x}\;\;\;\text{on $\partial\Omega$}.
    \end{equation}
    Next, substitute the expressions for $\bm{v}_n$, $\bm{\g}_n$, and $\bm{\o}_n$ from (\ref{eq:scaledQties-2}) into (\ref{eq:wkForm_coupledSyst-1}) with $\d_n F$ in place of $f$. Taking, in particular, arbitrary $\bm{\f}\in\C(\Omega)$ and using properties (\ref{eq:props_of_components}) we get
    \begin{equation}\label{eq:initial_est_wAn}
        2\nu(\bm{u}_n,\bm{\f})_{\H(\Omega)}=\d_n A_n+\overline{F\bm{b}_n}\cdot\bm{\g}_{\bm{\f}}+\bm{r}\times\overline{F\bm{b}_n}\cdot\bm{\o}_{\bm{\f}},
    \end{equation}
where
\begin{equation}\label{eq:AnDef}
\begin{aligned}
    &A_n:=\left((\bm{u}_n-\bm{\xi}_n-\bm{\zeta}_n\times\bm{x})\cdot\nabla\bm{u}_n,\bm{\f}\right)_{L^2(\Omega)}+(\overline{(\w_n-\bm{\chi}_n-\bm{\eta}_n\times\bm{x})\cdot\nabla\w_n},\bm{\f})_{L^2(\Omega)} \\
    &\;\;\;\;\;\;\;\;\;\;\;\;\;\;\;\;\;\;\;\;\;\;\;\;\;\,+M\left( \bm{\zeta}_n\times\bm{\xi}_n+\overline{\bm{\eta}_n\times\bm{\chi}_n}\right)\cdot{\bm{\g}_{\bm{\f}}}+\left[\bm{\zeta}_n\times(\textbf{I}\cdot\bm{\zeta}_n)+\overline{\bm{\eta}_n\times(\textbf{I}\cdot\bm{\eta}_n)}\right]\cdot\bm{\o}_{\bm{\f}}.
\end{aligned}
\end{equation}
Then, comparing (\ref{eq:initial_est_wAn}) with (\ref{eq:wkForm_lnrProb}), one immediately sees from (\ref{eq:propConvergences}$)_1$ and (\ref{eq:PropConvgces}$)_1$ that the lemma is proved once we show that $A_n$ is bounded uniformly in $n$; indeed, then we can pass to the limit in (\ref{eq:initial_est_wAn}) as $n\rightarrow\infty$ and, thanks to (\ref{eq:trace-result-in-prop-thm}), use the uniqueness property of Lemma \ref{thm:eom_steady_lowRe2}. To that end, we need to deduce some uniform estimates.
First, from (\ref{eq:energyForWRFinal_coupledSyst}), and (\ref{eq:est_on_un}) it follows that
\begin{equation}\label{eq:est_for_chin-etan}
    \|\bm{\chi}_n\|_{L^2(0,T)} + \|\bm{\eta}_n\|_{L^2(0,T)} \leq \kappa.
\end{equation}
Similarly, choosing $R>\text{diam}\,\mathcal{B}$ sufficiently large so that $\Omega_R\supset\text{supp}\,\nabla\bm{\f}$, again from (\ref{eq:energyForWRFinal_coupledSyst}) and (\ref{eq:est_on_un}) with the help of the Sobolev and H\"older inequalities, we deduce
    \begin{equation}\label{eq:estimate-for-wn}
            \|\w_n\|_{L^2(0,T;L^2(\Omega_R))} \leq \frac{1}{\d_n}\|\bm{v}_n\|_{L^2(0,T;L^2(\Omega_R))}+\|\bm{u}_n\|_{L^2(0,T;L^2(\Omega_R))} \leq c_1\kappa,
    \end{equation}
where $c_1=c_1(R)>0$. Then, employing in (\ref{eq:AnDef}) the uniform bounds (\ref{eq:est_on_un}),  (\ref{eq:est_for_chin-etan}), and (\ref{eq:estimate-for-wn}) in combination with Lemma \ref{thm:KornsIdentity_coupledSyst} and H\"older inequality, we easily prove that $A_n$ is indeed bounded, thus completing the proof of the lemma.
\end{proof}

With the help of Lemma \ref{thm:thrustApprox}, we are in a position  to prove the main result of this section. \\

\begin{thm}\label{thm:suffCond4prop}
    Let $(\bm{v},\bm{\g},\bm{\o},\bm{b})$ be a weak solution to problem (\ref{eq:eom_coupled_system}) corresponding to the force $f\in L_\text{per}^{\infty}(\mathbb{R})$, where $\overline{f}=:\delta\, \overline{F}\neq 0$. Then, if 
\begin{equation}\label{SP} \widehat{\bm{b}}\neq\textbf{C}\cdot\bm{\Uptheta}^{-1}\cdot(\bm{r}\times\widehat{\bm{b}})\,,\end{equation} 
necessarily $\overline{\bm{\g}}\neq \textbf{0}$; that is, $\mathcal{B}$ experiences propulsion. Precisely, there is $\delta_0>0$ such that
    \begin{equation}\label{*}
        \overline{\bm{\g}} = \delta\overline{F}\left(\textbf{K}-\textbf{C}\cdot\bm{\Uptheta}^{-1}\cdot\textbf{C}^{\top}\right)^{-1}\cdot\left(\widehat{\bm{b}}-\textbf{C}\cdot\bm{\Uptheta}^{-1}\cdot(\bm{r}\times\widehat{\bm{b}})\right) +{\sf R}(\delta)\,,\ \ \mbox{for all $\delta\in(0,\delta_0)$}\,,
\end{equation}
where
\begin{equation}\label{**}
\lim_{\delta\to 0}\frac1\delta\,{\sf R}(\delta)=0\,.
\end{equation}
\end{thm}
\begin{proof} For $\bm{a}\in \mathcal H(\Omega)$, we set
$$
\bfcalf(\bm{a}):=\sum_{i=1}^3(\bm{a},\bm{h}^{(i)})_{\mathcal{H}(\Omega)}\bm{e}_i\,,\ \ \bfcalg(\bm{a}):=\sum_{i=1}^3(\bm{a},\bm{H}^{(i)})_{\mathcal{H}(\Omega)}\bm{e}_i\,. 
$$
Dot-multiplying both sides of (\ref{eq:eom_aux-1})$_1$ by $\bm{u}_0$, integrating by parts over $\Omega$ and taking into account (\ref{eq:eom_avged-Lnr-coupled_system})$_3$ and (\ref{matrices}), we get
\begin{equation}\label{26}
    \bfcalf(\bm{u}_0)=\textbf{K}\cdot\bm{\xi}_0+\textbf{C}\cdot\bm{\zeta}_0.
\end{equation} 
Likewise, by dot-multiplying this time both sides of (\ref{eq:eom_aux-2})$_1$ by $\bm{u}_0$, integrating by parts over $\Omega$ and using again  (\ref{eq:eom_avged-Lnr-coupled_system})$_3$ and (\ref{matrices}), it follows that
\begin{equation}\label{27}
    \bfcalg(\bm{u}_0)=\textbf{C}^{\top}\cdot\bm{\xi}_0+\bm{\Uptheta}\cdot\bm{\zeta}_0.
\end{equation} 
Repeating the above procedure with $\bm{u}_\delta$ (defined in Lemma \ref{thm:thrustApprox}) in place of $\bm{u}_0$, and recalling (\ref{26})--(\ref{27}), we thus deduce
\begin{equation*}
\begin{aligned}
    \bfcalf(\bm{u}_\delta-\bm{u}_0)&=\textbf{K}\cdot(\bm{\xi}_\delta-\bm{\xi}_0)+\textbf{C}\cdot(\bm{\zeta}_\delta-\bm{\zeta}_0)\\
    \bfcalg(\bm{u}_\delta-\bm{u}_0)&=\textbf{C}^{\top}\cdot(\bm{\xi}_\delta-\bm{\xi}_0)+\bm{\Uptheta}\cdot(\bm{\zeta}_\delta-\bm{\zeta}_0),
\end{aligned}
\end{equation*}
which, in turn, furnishes
$$
\bm{\xi}_\delta=\bm{\xi}_0+\left(\textbf{K}-\textbf{C}\cdot\bm{\Uptheta}^{-1}\cdot\textbf{C}^{\top}\right)^{-1}\cdot\left(\bfcalf(\bm{u}_\delta-\bm{u}_0)-\textbf{C}\cdot\bm{\Uptheta}^{-1}\cdot\bfcalg(\bm{u}_\delta-\bm{u}_0))\right)\,.
$$
Therefore, (\ref{*})--(\ref{**}) follows from this equation by taking into account (\ref{def:wkSoltoStokesProb-1})$_1$ and (\ref{eq:PropConvgces})$_2$.
\end{proof}
\medskip
\begin{remark}
In case (\ref{SP}) is violated, it may happen that propulsion takes place at an order in $\delta$ higher than 1. This possibility is investigated in \cite{MMK} at the order of $\delta^2$.
\end{remark}
\begin{remark}\label{rem}
It is interesting to consider the counterpart of Theorem \ref{thm:suffCond4prop} in the case when $\mathcal B$ is constrained to execute a translational motion only. This can be achieved by applying a suitable torque on $\mathcal B$ to prevent rotational motion. We will only sketch the analysis, referring to \cite{MMK} for full details.  
In such a case, equations (\ref{eq:eom_coupled_system}) reduce to
\begin{equation}\label{eq:coupled_system}
\begin{aligned}
    \left.\begin{array}{c}
        \displaystyle
        \frac{\partial\bm{v}}{\partial t}+(\bm{v}-\bm{\gamma})\cdot\nabla\bm{v}=\div\,\textbf{T}(\bm{v},p) \\
        \;\;\;\;\;\;\;\;\;\;\;\;\;\;\;\;\;\;\;\;\;\div\,\bm{v}=0 \\
    \end{array}\right\}&\;\;\;\;\;\;\;{\text{in}\;\Omega\times\mathbb{R}} \\
    \bm{v}=\bm{\gamma}&\;\;\;\;\;\;\;\text{on}\;\partial\Omega\times\mathbb{R} \\
    \lim_{|\bm{x}|\rightarrow\infty}\bm{v}(\bm{x},t)=\textbf{0}\;&\;\;\;\;\;\;\;\text{in}\;\mathbb{R}. \\
                \displaystyle
        M\dot{\bm{\g}}=f\widehat{\bm{b}}-\int_{\partial\Omega}\textbf{T}(\bm{v},p)\cdot\bm{n}\,\text{d}S&\;\;\;\;\;\;\;{\text{in}\;\mathbb{R}}\,.
\end{aligned}
\end{equation}
For this problem, one can show a result of existence of $T$-periodic weak solutions $(\bm{v},\bm{\gamma})$ --in the sense of Definition \ref{def:wkSol-1}-- entirely analogous to Theorem \ref{thm:exist_wk_sol_coupledSyst}, with corresponding  estimate
$$
    \|\bm{v}\|_{L^2(0,T;\H(\Omega))}+\|\bm{\gamma}\|_{L^2(0,T)}\leq C\|f\|_{L^2(0,T)}\,.
$$
Employing the latter in combination with the scaling argument presented above, one can prove that the rescaled averaged solution $(\bm{\mathsf{u}}_\delta:=\frac1\delta\overline{\bm{v}},\bm{\mu}_\delta:=\frac1\delta\overline{\bm{\gamma}})$  converges, as $\delta\to 0$, to the unique solution to the problem
\begin{equation}\label{30}
\begin{aligned}
    \left.\begin{array}{c}
        \displaystyle \div\textbf{T}(\bm{\mathsf{u}}_0,{\sf p}_0)=\bm{0} \\
        \;\,\;\;\;\;\;\;\;\;\div\,\bm{\mathsf{u}}_0=0 \\
    \end{array}\right\}&\;\;\;\;\;\;\;{\text{in}\;\Omega} \\
    \bm{\mathsf{u}}_0=\bm{\mu}_0&\;\;\;\;\;\;\;\text{on}\;\partial\Omega \\
    \int_{\partial\Omega} \textbf{T}(\bm{\mathsf{u}}_0,{\sf p}_0)\cdot\bm{n}\;\text{d}S=\overline{F}\widehat{\bm{b}}\,.&
\end{aligned}
\end{equation}
If we dot-multiply both sides of (\ref{30})$_1$ by $\bm{h}^{(i)}$, integrate by parts over $\Omega$ and take into account (\ref{30})$_4$, we get
$$
\overline{F}\,\widehat{\bm{b}}=(\bm{\mathsf{u}}_0,\bm{h}^{(i)})_{\mathcal{H}(\Omega)}\,,\ \ i=1,2,3\,.
$$
Similarly, dot-multiplying both sides of (\ref{eq:eom_aux-1})$_1$ by $\bm{u}_0$, integrating by parts over $\Omega$ and using (\ref{30})$_3$ and (\ref{matrices})$_1$, we infer
$$
\big(\bm{\mu}_0\cdot\textbf{K}\big)_i=(\bm{\mathsf{u}}_0,\bm{h}^{(i)})_{\mathcal{H}(\Omega)}\,,\ \ i=1,2,3\,.
$$
Thus, combining the last two displayed equations we conclude
$$
\bm{\mu}_0=\overline{F}\,\textbf{K}^{-1}\cdot\widehat{\bm{b}}
$$
Thus, adapting to our case the proof of Theorem \ref{thm:suffCond4prop}, one can show that
$$
\overline{\bm{\gamma}}=\delta\,\overline{F}\,\textbf{K}^{-1}\cdot\widehat{\bm{b}}+o(\delta)\,\ \mbox{as $\delta\to0$}\,,
$$
which furnishes that, already at the first order in $\delta$, it is $\bm{\g}\neq\bm{0}$ provided only $\overline{F}\neq 0$, no matter the point where $\bm{\mathsf{f}}$ is applied and shape or physical properties of body $\mathcal B$. 
\end{remark}

\hfill

\section{An Example for the Sphere}\label{ssec:example_sphere}
In this section, we show that a homogeneous sphere (that is, of uniformly distributed mass) always experiences propulsion under the action of a time-periodic force $\bm{f}:=f\widehat{\bm{b}}$ if $\overline{f}\neq 0$, regardless of the location at which this force is applied. However, given a \textit{non-homogeneous} sphere with center of mass $G$, chosen suitably different than its geometric center, we show that there is at least one point $P$ such that the force $\bm{f}$, applied at $P$, does \textit{not} induce propulsion at the order of $\d$, even if $\overline{f}\neq 0$.

First, let $\mathcal{S}\subset\R^3$ be a homogeneous sphere with radius $a$ and geometric center $R$. Then, $R$ coincides with the center of mass $G$ which, by the convention outlined in Section \ref{ssec:intro}, is also taken to be origin of our frame-of-reference $\mathcal{F}:=\{G;\bm{e}_1,\bm{e}_2,\bm{e}_3\}$. It is known that, for such a sphere, the matrices $\textbf{K}, \bm{\Uptheta}$, and $\textbf{C}$, introduced in Lemma \ref{thm:eom_aux}, can be taken (upon possible rotation of $\mathcal{F}$) to be as follows (see equations (5-2.22), (5-3.13), and \textit{Case 3} in \cite[Section 5-5]{happel1983low}):
\begin{equation}\label{eq:formOfK-Theta-C-Sphere-homog}
    \textbf{K}=6\pi a\textbf{1},\;\;\;\;\;\bm{\Uptheta}=8\pi a^3\textbf{1},\;\;\;\;\;\text{and}\;\;\;\;\;\textbf{C}=\textbf{O}\,.
\end{equation}
Then, substituting (\ref{eq:formOfK-Theta-C-Sphere-homog}$)_3$ into the relation
\begin{equation}\label{eq:condPos4prop}
    \widehat{\bm{b}}=\textbf{C}\cdot\bm{\Uptheta}^{-1}\cdot(\bm{r}\times\widehat{\bm{b}}),
\end{equation}
gives $\widehat{\bm{b}}=\textbf{0}$, which is certainly not true. Hence, by Theorem \ref{thm:suffCond4prop}, if $\overline{f}\neq 0$, then $\mathcal{S}$ always propels.

Now, let us construct a non-homogeneous sphere $\mathcal{S}'\subset\R^3$ as follows: modify the mass distribution of $\mathcal{S}$ such that its center of mass is now located at a point $G'$ along the axis $\bm{e}_1$ at a distance $d\in(0,a)$ from $R$, noting that this point also lies on $\bm{e}_1$ (see Figure \ref{spherefig}). In this case, from relations (5-4.10) and (5-4.12) of \cite{happel1983low}, we have

    \InsertBoxR{0}
    {
    \begin{minipage}{0.29\textwidth}\centering
        \includegraphics[width=1\textwidth]{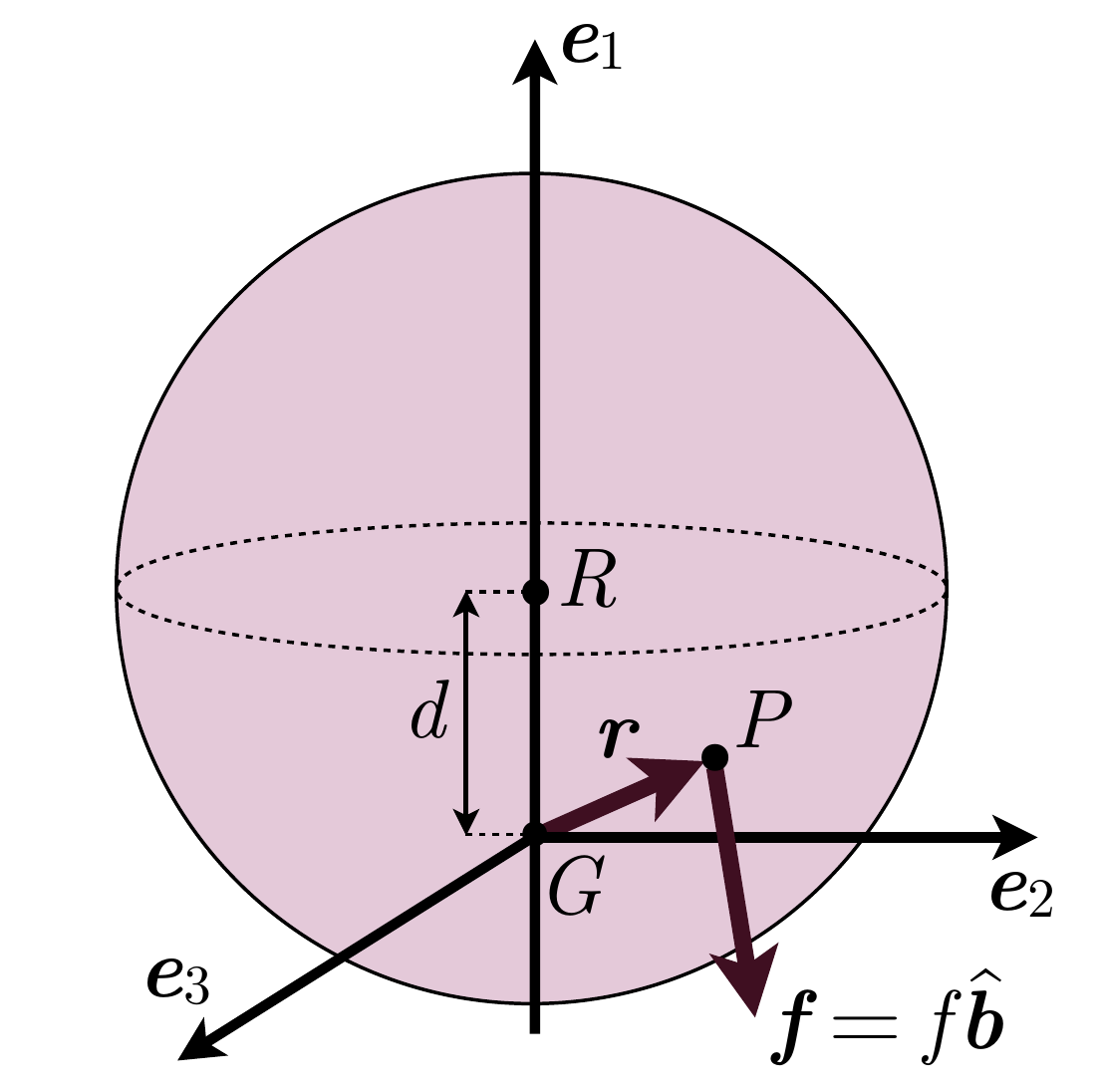}
    \captionof{figure}{\centering Schematic of Sphere $\mathcal{S}'$}
        \label{spherefig}
    \end{minipage}
    }
 
\begin{equation}\label{eq:formOfTheta-C-Sphere-nonhomog}
\begin{aligned}
    \bm{\Uptheta}&=2\pi a\begin{pmatrix}
        4a^2 & 0 & 0 \\
        0 & 4a^2+3d^3 & 0 \\
        0 & 0 & 4a^2+3d^3
    \end{pmatrix} \\
\textbf{C}&=6\pi a\begin{pmatrix}
            0 & 0 & 0 \\
            0 & 0 & -d \\
            0 & d & 0
        \end{pmatrix}.
\end{aligned}
\end{equation}
Then, one easily finds that (\ref{eq:condPos4prop}) is verified with (\ref{eq:formOfTheta-C-Sphere-nonhomog}) along with
$$
    \widehat{\bm{b}}:=(0,\widehat{b}_2,\widehat{b}_2),
$$
for any $\widehat{b}_2,\widehat{b}_2\in\mathbb{R}$ (not both zero) and
$$
    \bm{r}:=\left(r_1, r_2, r_3\right),
$$
with
$$
    r_1:=-\frac{1}{12\pi^2 a^2d(4a^2+3d^3)}\,,\  \ \ |r_i|<a-d \ \ \mbox{for all $i=1,2,3$}\,, 
$$
where, of course, these vectors are expressed in the reference frame $\mathcal{F}':=\{G';\bm{e}_1,\bm{e}_2,\bm{e}_3\}$. Then by Theorem \ref{thm:suffCond4prop}, at the order $\d$, the sphere $\mathcal{S}'$ does not propel \textit{regardless} of the value of $\overline{f}$. In other words, $\mathcal{S}'$ exhibits a purely oscillatory motion. 

\hfill

\section*{Acknowledgments}
Work partially supported by National Science Foundation (US) Grant DMS-2307811

\hfill

\bibliographystyle{ieeetr} 
\bibliography{References.bib}

\end{document}